\documentclass[12pt]{amsart}

\usepackage{amsmath,amssymb,latexsym} 
\usepackage{fullpage}
\usepackage{etoolbox}
\usepackage{enumitem}
\usepackage{color}
\usepackage[colorinlistoftodos,bordercolor=orange,backgroundcolor=orange!20,linecolor=orange,textsize=scriptsize]{todonotes}
\usepackage{enumitem}

\newtheorem{theorem}{Theorem} 

\newtheorem{lemma}{Lemma}

\theoremstyle{definition}

\theoremstyle{remark}

\newtheorem*{proof-claim}{Proof}

\newenvironment{changemargin}[2]{\begin{list}{}{%
\setlength{\topsep}{0pt}%
\setlength{\leftmargin}{0pt}%
\setlength{\rightmargin}{0pt}%
\setlength{\listparindent}{\parindent}%
\setlength{\itemindent}{\parindent}%
\setlength{\parsep}{0pt plus 1pt}%
\addtolength{\leftmargin}{#1}%
\addtolength{\rightmargin}{#2}%
}\item }{\end{list}}

\AtBeginEnvironment{proof-claim}{\vspace{-0.5cm}\footnotesize\begin{changemargin}{0.5cm}{0.5cm}}
\AtEndEnvironment{proof-claim}{\end{changemargin}}

\def\L{\mathsf{L}}
\def\R{\mathbb{R}}
\def\Z{\mathbb{Z}}
\def\Fam{\mathcal{F}}
\def\L{\mathcal{L}}
\def\S{\mathbb{S}}
\def\C{\mathbb{C}}

 \title{Complex homothetic sections and projections\\ through a Helly type Theorem for cosets of $\S^1$}

\author{Jorge Luis Arocha}
\author{Javier Bracho }
\author{Luis Montejano }

\address{UNAM, Mexico}
\email{luis@matem.unam.mx}
\date{\today}

\begin{document}
\maketitle
\begin{abstract}
We prove that two closed subsets of complex space $\C^n$ with  corresponding complex homothetic sections (projections)
are complex homothetic. The proof uses a new Helly-type theorem for cosets of closed subgroups of  $\S ^1$.   
\end{abstract}

\section{Introduction }

{\em Let $V$ be an $n$-dimensional vector space  and let $K_1, K_2$ be two closed sets, $n\geq3.$
Suppose that for every 2-dimensional subspace (hyperplane) $\Lambda$ the corresponding sections 
$$ \Lambda\cap K_1  \mbox{  and }  \Lambda\cap K_2$$
are ``equivalent''.   Is it true that 
$$K_1  \mbox{  is ``equivalent'' to }  K_2?$$
Where the meaning of ``equivalent'' will be clarified in the sequel. }

\medskip

Suppose for example that $K_1, K_2$ are convex bodies containing the origin  in euclidean $n$-space, $n\geq 3$,  and the meaning of ``equivalent'' is congruence. 
That is,  suppose that for every hyperplane $\Lambda$ through the origin  there is a linear isometry  $g\in O(\Lambda)$ and a vector $v\in \Lambda$ (depending on $\Lambda$) such that 
$$g(\Lambda\cap K_1) + v = \Lambda\cap K_2.$$
Is it true that $K_1$ is either a translated copy of $K_2$ or $-K_2$?

V.P. Golubyatnikov \cite{ G} and D. Ryagobin  \cite{R0}, \cite{R1}, \cite{R3}, \cite{R2}, have studied this problem extensively.  Perhaps the most representative result in this direction is   Ryagobin's Theorem in \cite{R0}.  Suppose  $K_1, K_2$ are convex bodies in euclidean $n$-space such that for every plane 
$\Lambda$ through the origin  there is a rotation  $r$ of $\Lambda$ about the origin  (depending on $\Lambda$) for which 
$r(\Lambda\cap K_1)  = (\Lambda\cap K_2)$,
then either $K_1=K_2$ or $K_1=-K_2$.

Ryagobin's proof rests heavily on the fact  that there are only a countable number of closed subgroups of $SO(2)$. The analogous is not true for $SO(k)$, $k>2$, and this is a reason why the result is unknown for $k$-subspaces in $\R^n$. 

In the charming paper \cite{R}, Rogers proved the following. Suppose $K_1, K_2$ are convex bodies in euclidean $n$-space, $n\geq 3$, such that for every hyperplane $\Lambda$ through the origin 
$$\Lambda\cap K_1  \mbox{  is positively homothetic to } \Lambda\cap K_2,$$
then $K_1$ is positively homothetic to $K_2.$

\smallskip
Note that the center of each homothecy between $\Lambda\cap K_1$ and $\Lambda\cap K_2$ is not  fixed and of course depends on the hyperplane $\Lambda$. This result was crucial in many characterizations of ellipsoids in the second half of the 20th century.
 
Studying complex ellipsoids in \cite{ABMCE}, it was necessary to prove a complex version of Rogers' Theorem for real homotheties. More precisely, suppose  $K_1, K_2$ are convex bodies in complex $n$-space,  $\C^n$, $n\geq 3$,  such that for every complex hyperplane $\Lambda$ through the origin 
$$\Lambda\cap K_1  \mbox{  is positively homothetic to } \Lambda\cap K_2,$$
then $K_1$ is positively homothetic to $K_2.$

It would be interesting to know if the following complex version of Ryagobin's Theorem is true.  
Suppose two convex bodies $K_1, K_2$  in complex $n$-space  $\C^n$ have the property that for every complex plane $\Lambda$, there is $g\in U(\Lambda)$  such that 
$$g(\Lambda\cap K_1)= \Lambda\cap K_2\,.$$
Is there $z\in \S^1$ with $|z|=1$ such that 
 $zK_1=K_2$?
 
 It is precisely in this direction that we obtain the following theorem as a special case of the above. We note that our result is true not only for convex subsets but also for closed subsets of $\C^n$ (including of course finite sets). 
  
 \begin{theorem}\label{mainsec}
 Let $K_1, K_2$ be two closed subsets of $\C^n$, $n\geq3$.
Suppose that for every complex plane $\Lambda$ through the origin there is a complex number $z_\Lambda\in\C$ such that 
$$ z_\Lambda(\Lambda\cap K_1)=\Lambda\cap K_2.$$
Then, there is $z\in \C$, such that 
$$zK_1  = K_2.$$
\end{theorem}
 
 The analogue of Theorem \ref{mainsec} in the real setting is essentially  Lemma 1 of \cite{R0}.  
 
  The proof of  Theorem \ref{mainsec} makes use of the following Helly-type theorem for cosets of closed subgroups of  $\S ^1$. 
  
  \begin{theorem} \label{Teo_Helly_orbital}
Let $\Fam$ be a pairwise intersecting family of cosets of closed subgroups of  $\S^1$.  Then $\Fam$  has non-empty intersection. 
  \end{theorem}
  
  It would be interesting to know which groups have this Helly property. That is, for which groups $G$ there exists an integer $m>2$ with the property that if $\Fam$ is a family of m+1 cosets of closed subgroups of $G$ and each $m$ members de $\Fam$ is intersecting, then $\Fam$ is intersecting.  Indeed, the integers $\Z$ have this property. We will prove that every three pairwise intersecting arithmetic sequences is intersecting (see Theorem~\ref{Teo_Helly_integers}). 
    
  Similarly, using Theorem \ref{Teo_Helly_orbital}, we prove the corresponding result for projections.  
  Let $\Lambda$ be a complex subspace of $\C^n$, we denote by $\pi_\Lambda:\C^n\to \Lambda$, the orthogonal projection.
  
   \begin{theorem}\label{mainproj}
 Let $K_1, K_2\subset \C^n$ be either two closed convex sets or two finite sets, $n\geq3$. 
Suppose that for every complex plane $\Lambda$ through the origin there is a complex number $z_\Lambda\in\C$ such that 
$$ z_\Lambda(\pi_\Lambda(K_1))=\pi_\Lambda(K_2).$$
Then, there is $z\in \C$, such that 
$$zK_1  =  K_2.$$
\end{theorem}

Regarding projections of convex bodies in $\R^5$  whose projections are $SU(2)$-congruent, Ryagobin obtained interesting results in \cite{R3}.

\section{The group of complex homotheties}

From now on, unless otherwise stated, all our lines, planes, hyperplanes pass through the origin.
 
Let $V$ be a complex vector space and let $K\subset V$ be a closed subset. Define the {\em group of  homotheties of } $K$ as:
$$H(K)=\{\,z\in\S^1\,|\, zK=K\,\},$$
which is a closed subgroup of $\S^1$.
For every subspace $\Delta\subset \C^n$, 
\begin{equation}
 H(K)\mbox{  is a subgroup of  } H(\Delta\cap K),
 \end{equation}
 \begin{equation}\label{Eq:Sub_proy}
 H(K)\mbox{  is a subgroup of  } H(\pi_\Delta (K)).
 \end{equation}
 
 Note that (1) is trivially true. On the other hand, if $\pi_\Delta:\C^n\to \Delta$ is the orthogonal projection, then 
 $z\pi_\Delta(x)=\pi_\Delta(zx)$, for every 
 $z\in \C$ and every $x\in \C^n$.  Hence, $zK=K$  implies $z\pi_\Delta K=\pi_\Delta K$ and consequently (2) is true.

Now, if $K_1, K_2$ are closed subsets of $V$,  let us define 
$$H(K_1, K_2)=\{z\in\S^1|zK_1=K_2\}.$$

For every $z_0\in H(K_1,K_2)$, \begin{equation}\label{Eq:Cosets}
z_0H(K_1)=H(K_1,K_2)\,.
\end{equation}
Therefore, if $H(K_1,K_2)$ is nonempty, it is a coset of the closed subgroup $H(K_1)$ of $\S^1$.
 
  Let $\L$ denote the set of complex lines through the origin in $\C^n$. Of course, it is the complex projective space of dimensión $n-1$, but we need to refer to its elements as complex lines.
  
    \begin{lemma}\label{lemlines}
  Let $K\subset \C^n$, $n\geq3$, be a closed subset.  Then
 $$H(K)=\bigcap_{L\in\L}H(L\cap K).$$
\end{lemma}
\proof  By (1), $H(K)\subset \bigcap_{L\in\L}H(L\cap K)$. 
 Suppose now, $z\in \bigcap_{L\in\L}H(L\cap K)$,  that is, $z(L\cap K)=L\cap K$, for every line $L\in\L$. Consequently, $zK=K$, which implies that $z\in  H(K)$ as we wished.
\qed 
  
\section{Helly type results}

For the proof of our main Helly type theorem (Theorem~\ref{Teo_Helly_orbital}), we need two lemmas.

\begin{lemma}\label{Helly-series}
Let $a, b \in \Z$, then $a \Z$, $b \Z$ and $a + (b-a)\Z$ have non empty intersection.
\end{lemma}

\proof Let $d$ be the greatest common divisor of $a$ and $b$. There exist $a^\prime, b^\prime, x, y \in \Z$ such that
$$a=d a^\prime \quad \text{,}\quad  b=d b^\prime$$
and 
$$x a^\prime + y b^\prime=1\,.$$
Adding and subtracting $x b^\prime$
$$x (a^\prime-b^\prime) + (x+y) b^\prime=1\,.$$
From which
$$(x+y) b^\prime=1+x (b^\prime-a^\prime)\,,$$
and multiplying by $d a^\prime=a$, we get
$$(x+y) d a^\prime b^\prime=a+x a^\prime (b-a)\,.$$
The left hand side is in $a \Z\cap b \Z$; the right hand side is in $a + (b-a)\Z$. \qed

\begin{lemma}\label{Helly-3} Let $A, B, C$ be cosets of closed subgroups of $\S^1$ such that each pair is intersecting, then they intersect.
\end{lemma}

\proof 
If one of the given cosets is $\S^1$, we are done, so we may assume that the three cosets are finite. 
Without loss of generality, we can suppose that $A$ and $B$ are subgroups of $\S^1$ (so that ${1\in A\cap B}$), for otherwise, take $z\in A\cap B$, and consider $z^{-1}A,  z^{-1}B, z^{-1}C$.

By hypothesis, there exist $\alpha\in A\cap C$ and $\beta\in B\cap C$. Observe that $\alpha\langle \beta\alpha^{-1} \rangle\subset C$, where we denote by $\langle z \rangle$ the subgroup of $\S^1$ generated by $z\in\S^1$.

Let $p=order(\alpha)$ and $q=order(\beta)$. So that $\langle \alpha \rangle\subset A$ and $\langle \beta \rangle\subset B$ are the cyclic groups of order $p$ and $q$  respectively.

Let $m$ be the minimum common multiple of $p$ and $q$. And let
\begin{align*}
\varphi: \Z&\to\S^1  \\
               \varphi(n)&=exp(n\frac{2\pi}{m})
\end{align*}
where $exp:\R\to\S^1$ is the standard exponential map. By construction, there exist $a, b\in\Z$ such that $\varphi(a)=\alpha$ and $\varphi(b)=\beta$.

By Lemma~\ref{Helly-series}, there exists $c\in a\Z\cap b\Z\cap (a+(b-a)\Z)$, so that
$$\varphi(c)\in\langle \alpha \rangle\cap\langle \beta \rangle\cap\alpha\langle \beta\alpha^{-1} \rangle \subset A\cap B\cap C\,.$$
\qed

\medskip
\noindent {\bf Proof of Theorem \ref{Teo_Helly_orbital}}.   Let $\Fam$ be a pairwise intersecting family of cosets of closed subgroups of $\S^1$. 

First, we prove the case where $\Fam$ is finite by induction on 
the cardinality, $r$, of $\Fam$. For $r=3$ it is Lemma~\ref{Helly-3}. Suppose the Theorem true for families of size $r-1$. Consider $A\in\Fam$ and let 
$\Fam^\prime=\{\,A\cap X\,|\,X\in\Fam\setminus\{A\}\,\}$. By Lemma~\ref{Helly-3}, every pair of elements of $\Fam^\prime$ intersect and by induction they all have non-empty intersection which is also the intersection of all the elements of $\Fam$. 

Finally, consider the general case. Suppose that $\Fam$ is not intersecting. Then the family of its complements cover $\S^1$. They are open sets, so that the compactness of $\S^1$ gives us a finite open subcover. By the preceding paragraph, in the corresponding finite family of closed cosets of groups of $\S^1$ there must be a non intersecting pair, which is a contradiction.
\qed

Again using Lemma \ref{Helly-series},  it is easy to prove that not only  the group $\S^1$ has this Helly property  but also the group $\Z$ for finite families.

 \begin{theorem} \label{Teo_Helly_integers}
  Let ${\Fam}$ be a finite pairwise intersecting family of cosets of subgroups of  $\Z$.  Then ${\Fam}$  has non empty intersection. 
  \end{theorem}
  
  \proof Following the proof of Theorem \ref{Teo_Helly_orbital}, it is enough to consider the case when ${\Fam}=\{A,B,C\}$ consists of three cosets of subgroups of $\Z$. Therefore $A, B$ and $C$ are arithmetic  sequences and without loss of generality $A=a\Z$, and $B=b\Z$.
Assume $a'\in C\cap A$ and $b'\in C\cap B$,  Consequently, $a'\Z\subset A$,   $b'\Z\subset B$ and 
 $a'+(b'-a')\Z \subset C$. By Lemma \ref{Helly-series}, we have that ${\Fam}=\{A,B,C\}$ is intersecting.  This concludes the proof of the theorem. \qed
 
  \section{ Complex homothetic sections and projections }
  
  Let $K_1, K_2$ be two non-homothetic complex ellipsoids centered at the origin in $\C^n$.
  Hence, for every line $L$, $L\cap K_1$ and $L\cap K_2$ are disks centered at the origin and consequently 
 there is $z_L$ such that 
 $$ z_L(L\cap K_1)= L\cap K_2.$$
 Nevertheless, $zK_1\not=K_2$, for  $z\in \C$.  This shows that Theorem 1 is false if we replace complex planes with complex lines. 
  
 For the proof of Theorem \ref{mainsec}, we will start by noting that in its hypothesis we can assume, without loss of generality, that the module of all $z_\Lambda$ is 1.
  
\begin{lemma}\label{lemwlg} 
Let  $K_1, K_2$ be closed subsets of $\C^n, n\geq3$ and
suppose that for every complex plane $\Lambda$   there is $z_\Lambda\in \C$, such that 
$$ z_\Lambda(\Lambda\cap K_1)=\Lambda\cap K_2.$$
Then, we may assume that
$ |z_\Lambda|$
is constant.
\end{lemma}

\proof Let us call a subset of $\C^n$ \emph{nontrivial} if it is neither empty nor the origin, and observe that for a nontrivial set $A\subset\C^n$, if $z_1A=z_2A$ then $|z_1|=|z_2|$. 
We may clearly assume $K_1$ is nontrivial, so that there exists a line $L$ such that $L\cap K_1$ is nontrivial.
If $\Lambda_1$ and $\Lambda_2$ are two complex planes which contain the complex line $L$,  then by hypothesis, 
$$z_{\Lambda_1}(L\cap K_1)=L\cap K_2=z_{\Lambda_2}(L\cap K_1)$$ 
and then, $|z_{\Lambda_1}|=|z_{\Lambda_2}|=r_L$.  
Now consider a plane $\Lambda$ that does not contain $L$. If $\Lambda\cap K_1$ is trivial, we may choose $z_{\Lambda}=r_L$. Otherwise, there is a line $L^\prime\subset\Lambda$ such that $L^\prime\cap K_1$ is nontrivial. Let  $\Lambda^\prime$ be the plane generated by $L$ and $L^\prime$. Then, 
$$z_{\Lambda}(L^\prime\cap K_1)=L^\prime\cap K_2=z_{\Lambda^\prime}(L^\prime\cap K_1)\,,$$ 
and therefore, $|z_{\Lambda}|=|z_{\Lambda^\prime}|=r_L$. \qed

\medskip

\noindent {\bf Proof of Theorem \ref{mainsec}}.  First, by Lemma \ref{lemwlg} we may assume that $|z_\Lambda|=1$ for every complex plane $\Lambda$. Let $L$ be a complex line and $L \subset \Lambda$, where $\Lambda$ is a complex plane. Then $ z_\Lambda(L\cap K_1)=(L\cap K_2).$ 
Therefore $z_\Lambda$ lies in the coset $  H(L\cap K_1,L\cap K_2)\subset \S^1$ which  is not empty.
Let 
$$\Fam=\{\,  H(L\cap K_1,L\cap K_2)\,\mid\, L\in\L \,\}\,.$$  
We shall prove that $\Fam$ is a pairwise intersecting family of cosets of closed subgroups of $\S^1$.  Let $L_1, L_2$ be two different complex lines and let $\Lambda$ be the complex plane generated by $L_1$ and $L_2$. Hence $$z_\Lambda\in  H(L_1\cap K_1,L_1\cap K_2)\cap H(L_2\cap K_1,L_2\cap K_2)\not=
\emptyset\,.$$ 
By Theorem~\ref{Teo_Helly_orbital} and  Lemma~\ref{lemlines}, there exists 
$$z\in \bigcap_{L\in\L}H(L\cap K_1,L\cap K_2)=H(K_1,K_2)\,.$$
Therefore, $zK_1=K_2$ as we wished. \qed

For the proof of Theorem \ref{mainproj} we need two lemmas.

\begin{lemma} \label{lemA}
 If $K_1, K_2\subset \C^n$ are either two closed convex sets or two finite sets such that  $\pi_L(K_1)=  \pi_L(K_2)$ for every complex line $L$, then $K_1=K_2$.
\end{lemma}

\proof  Suppose first $K_1, K_2\subset \C^n$ are two closed convex sets.  If $K_1\not=K_2$, then without loss of generality there is a point $x\in K_1$ and a real affine hyperplane $\Gamma$ such that $x$ is contained in the open semispace generated by $\Gamma$ and $K_2$ is contained in its complement. Let $\Delta$ be the unique complex affine hyperplane contained in $\Gamma$ and let $L$ be the complex line through the origin orthogonal to $\Delta$.  Consider $\pi_L:\C^n\to L$ be the orthogonal projection. Then $\pi_L(x)\in \pi_L(K_1)$ is contained in the open half plane generated by the real line $\pi_L(\Gamma)$ of $L$ and  $\pi_L(K_2)$ is contained in the complement, contradicting the assumption that $\pi_L(K_1)=\pi_L(K_2)$.  

Suppose now $K_1, K_2\subset \C^n$ are two finite sets. If they are different, then without loss of generality we may chose $x\in K_1\setminus K_2$. Since 
$K_1\cup K_2$ is finite there exists a complex line through the origin $L$ such that 
 $\pi_L\mid_{K_1\cup K_2}$ is injective. Then $\pi_L(x)\notin K_2$, contradicting the assumption that $\pi_L(K_1)=\pi_L(K_2)$.  
\qed 
  
 \begin{lemma}\label{lemI}
  Let $K, K_1, K_2 \subset \C^n$, $n\geq3$, be either closed convex sets or finite sets. Then
  \begin{enumerate}
 \item $H(K)=\bigcap_{L\in\L}H(\pi_L(K))\,,$
 \item $H(K_1, K_2)=\bigcap_{L\in\L}H(\pi_L(K_1), \pi_L(K_2)  ).$
\end{enumerate}
\end{lemma}

\proof  We know that $H(K)\subset \bigcap_{L\in\L}H(\pi_L(K))$  from  \eqref{Eq:Sub_proy}. 
 Suppose now, ${z\in \bigcap_{L\in\L}H(\pi_L(K))}$,  that is, $z\pi_L(K)=\pi_L(K)$, for every line $L\in\L$. 
 Consequently,  $\pi_L(zK)=\pi_L(K)$, for every line $L$ and therefore, by Lemma  \ref{lemA}, $zK=K$, which implies that $z\in  H(K)$ as we wished.  
 
Moreover,  for $z_0\in H(K_1, K_2)$,  by \eqref{Eq:Cosets} and item {\it (1)} of this lemma, we have:  
\begin{align*}
H(K_1, K_2)&=z_0H(K_1)=z_0\bigcap_{L\in\L}H(\pi_L(K_1))=\\
                    &=\bigcap_{L\in\L}z_0H(\pi_L(K_1))=\bigcap_{L\in\L}H(\pi_L(K_1), \pi_L(K_2)  ).
\end{align*}
\qed 

\medskip

\noindent {\bf Proof of Theorem \ref{mainproj}}. As in Lemma \ref{lemwlg}, let us first prove that we may assume that $|z_\Lambda|$
is constant and does not depend on $\Lambda$.  Remember from Section 2 that if $A$ is a closed subset of $\C^n$, $\Delta$ is a subspace of $\C^n$,  $z\in \C$ and $zA=A$, then $z\pi_\Delta(A)=\pi_\Delta(A)$.  

Now, call a line $L$ \emph{nontrivial} if $\pi_L(K_1)$ is nontrivial.
If  $L$ is a nontrivial complex line contained in the complex planes $\Lambda_1$ and  $\Lambda_2$, then $z_{\Lambda_1}(\pi_L(K_1))=\pi_L(K_2)=z_{\Lambda_2}(\pi_L(K_1))$ and so, $|z_{\Lambda_1}|=|z_{\Lambda_2}|.$  If $\Lambda_1$ and $\Lambda_2$ are two complex planes with nontrivial $\pi_{\Lambda_i}(K_1)$,  there is a third complex plane $\Lambda_3$ such that $\Lambda_i\cap\Lambda_3$ is a nontrivial complex line, $i=1,2$, so we obtain that $|z_{\Lambda_1}|= |z_{\Lambda_3}|=|z_{\Lambda_2}|$. 
Consequently, we may assume without loss of generality that $|z_\Lambda|=1$, for every complex plane.  

This time, let 
$$\Fam=\{  H(\pi_L(K_1), \pi_L (K_2))\mid L\in\L\}\,.$$  
We shall prove that $\Fam$ is a pairwise intersecting family of cosets of closed subgroups of $\S^1$.  Let $L_1, L_2$ be two different complex lines and let $\Lambda$ be the complex plane generated by $L_1$ and $L_2$. Hence, 
$$z_\Lambda\in H(\pi_{L_1}(K_1), \pi_{L_1} (K_2))\cap  H(\pi_{L_2}(K_1), \pi_{L_2} (K_2))\,.$$ 
By Theorem \ref{Teo_Helly_orbital}, there is 
$$z\in \bigcap_{L\in\L}H(\pi_L( K_1), \pi_L(K_2))\,.$$ 
Therefore,  by Lemma \ref{lemI}\it{(2)}, $z\in H(K_1,K_2)$. \qed

\medskip

\noindent {\bf Acknowledgments.} Luis Montejano acknowledges  support  from CONACyT under 
project 166306 and  from PAPIIT-UNAM under project IN112614. Javier Bracho acknowledges  support  from PAPIIT-UNAM under project IN109023.

\end{document}